\documentclass[a4papersize]{article}
\usepackage{amsmath,amsthm, amsxtra,amssymb,latexsym, amscd}
\usepackage[mathscr]{eucal}
\newtheorem{theorem}{Theorem}[section]
\newtheorem{corollary}[theorem]{Corollary}
\newtheorem{lemma}[theorem]{Lemma}
\newtheorem{proposition}[theorem]{Proposition}
\newtheorem{example}[theorem]{Example}

\newtheorem{remark}[theorem]{Remark}

\DeclareMathOperator{\depth}{depth}

\DeclareMathOperator{\Hom}{Hom}
\DeclareMathOperator{\Ann}{Ann}

\DeclareMathOperator{\Ass}{Ass}

\DeclareMathOperator{\Var}{Var}

\DeclareMathOperator{\Rad}{Rad}
\DeclareMathOperator{\Att}{Att}

\DeclareMathOperator{\Spec}{Spec}
\DeclareMathOperator{\p}{\frak p}

\DeclareMathOperator{\Ext}{Ext}
\begin{document}
\large
\centerline{\Large {\bf ASCENT AND DESCENT OF ARTINIAN MODULE STRUCTURES}}
\smallskip

  \centerline{\Large {\bf  UNDER FLAT BASE CHANGES}}
\medskip

\vskip 0.7cm
\centerline { TRAN DO MINH CHAU}
\centerline {  Thai Nguyen University of Education}
\centerline{Thai Nguyen, Vietnam}
\centerline {E-mail: chautdm@tnue.edu.vn}
\vskip 0.4cm

\centerline {LE THANH NHAN}
\centerline { Ministry of Education and Training, Hanoi, Vietnam}
\centerline {and Thai Nguyen University of Sciences, Thai Nguyen, Vietnam}
\centerline {E-mail: nhanlt2014@gmail.com}

\vskip 0.4cm
\smallskip
\vskip 1cm

 \noindent{\bf Abstract}  {\footnote{ {\it{Key words and phrases: }} Ascent and descent of module structure; Flat extension; Artinian module; Local cohomology module. \hfill\break {\it{2021 Subject  Classification: }} 13B40, 13E10, 13E05, 13D45. \hfill\break The work is supported by the Vietnam Ministry of Education and Training under grant number  B2022-TNA-25.}}.  Let $\varphi: R\rightarrow S$ be a flat local homomorphism between commutative Noetherian local rings. In this paper, the ascent and descent of Artinian module structures between $R$ and $S$ are investigated. For an Artinian $R$-module $A$, the structure of $A\otimes_RS$ is described. As an application, the Artinianess of certain local cohomology modules is clarified. 
 
\section{Introduction}  

\ \ \ Throughout this paper, let $\varphi: (R, \frak m)\rightarrow (S, \frak n)$ be a flat local homomorphism between commutative Noetherian local rings. The ascent and descent of finitely generated modules between $R$ and $S$ have attracted the interest of mathematicians, see \cite{AW}, \cite{FWW}, \cite{W}. The motivation for their investigations is the following vanishing result by A. J. Frankild and S. Sather-Wagstaff \cite{FW}: A finitely generated $R$-module $M$ is $\frak m$-adically complete if and only if $\Ext^i_R(\widehat R, M)=0$ for all $i\geq 1,$   if and only if $\Ext^i_R(F, M)=0$ for all flat $R$-modules $F$ and all integers $i\geq 1$. 

Assume that the induced map $R/\frak m\rightarrow S/\frak mS$ of $\varphi$ is an isomorphism, i.e. $\ell_R(S/\frak mS)=1$ (we note that the natural maps from $R$ to its $\frak m$-adic completion $\widehat R$ and to its Henselization $R^h$ satisfy this assumption). It is well known that if $M$ is a finitely generated $R$-module, then the  following five conditions are equivalent:

(a) $M$ has a $S$-module structure (with scalar multiplication $(s,x)\mapsto s\circ x$) that is compatible with its original $R$-module structure (i.e., $rx=\varphi(r)\circ x$ for all $r\in R, x\in M$);

(b) The natural map $M\to M\otimes_R S$ (sending $x$ to $x\otimes 1)$ is bijective;

(c) The natural map $\Hom_R(S, M)\to M$ (sending $f$ to $f(1))$ is bijective;

(d) $M\otimes_R S$ is a finitely generated as an $R$-module;

(e) $\Ext^i_R(S, M)=0$ for all $i>0.$

\noindent Moreover, if $N$ is a  finitely generated $S$-module, then $N$ is finitely generated as an $R$-module if and only if the natural map $N\to N\otimes_R S$ (sending $y$ to $y\otimes 1)$ is bijective, if and only if $N\otimes_R S$ is finitely generated as an $R$-module (see \cite{AW}, \cite{FWW}).

The purpose  of this paper is to investigate the ascent and descent of Artinian module structures between $R$ and $S$. We will show that in the case where  $\ell_R(S/\frak mS)=1$, the ascent and descent of Artinian module structures between $R$ and $S$ appear naturally, so that they are much different from the above ascent and descent of finitely generated module structures. 

\begin{proposition} Assume that $\ell_R(S/\frak mS)=1$. Let $A$ be an Artinian $R$-module, let $B$ be an Artinian $S$-module. Then the following statements are  true:

(a) $A$ has an Artinian $S$-module structure compatible with its original $R$-module structure, the natural maps $A\rightarrow A\otimes_R S$ and $\Hom_R(S, A)\rightarrow A$ are bijective, $A\otimes_R S$ is Artinian as an $R$-module  and $\Ext^i_R(S, A)=0$ for all integers $i>0$.  

(b) $B$ is an Artinian $R$-module by means of $\varphi$.  If we regard this $R$-module as a $S$-module mentioned in (a), then  we recover the original $S$-module structure on $B$.
\end{proposition}

As the main results, we first have the following characterizations for the descent of Artinian module structure in general case. Denote by $E(S/\frak n)$ the injective hull of the residue field $S/\frak n$ over $S$. Let $H^r_{\frak n}(S/\frak mS)$ be the top local cohomology module of $S/\frak mS$ with respect to $\frak n$, where $r=\dim(S/\frak mS)$. Note that $E(S/\frak n)$ and $H^r_{\frak n}(S/\frak mS)$ are Artinian $S$-modules. 

\begin{theorem}  The following statements are equivalent:

(a)  Each Artinian $S$-module is  Artinian as an $R$-module (by means of $\varphi$).

(b) $E(S/\frak n)$ is Artinian as an $R$-module.

(c) $H^r_{\frak n}(S/\frak mS)$ is Artinian as an $R$-module, where $r=\dim (S/\frak mS)$.

(d) $\ell_R(S/\frak m S)<\infty$.
\end{theorem}

In the statements (b), (c) of Theorem 1.2, the Artinian $S$-modules $E(S/\frak n)$ and $H^r_{\frak n}(S/\frak mS)$  can not be replaced by an arbitrary Artinian $S$-module (Example \ref{E:1}).

Secondly, we clarify the structure of $A\otimes_RS$, where $A\neq 0$ is an Artinian $R$-module. We note that $A\otimes_RS$ is an Artinian $S$-module if and only if $\dim (S/\frak mS)=0$ (Proposition \ref{P:1a}). If $\ell_R(S/\frak mS)<\infty$, we have the following ascent of Artinian module structures.

\begin{theorem} Assume that $\ell_R(S/\frak mS)=m$. Let $A$ be an Artinian $R$-module. Then  $A^m$ has a structure of Artinian $S$-module and there is an isomorphism $A\otimes_RS\cong A^m$ of $S$-modules. If we regard this $S$-module as an $R$-module by means of $\varphi$, then we recover the usual $R$-module structure on $A^m$.  
\end{theorem}

Theorem 1.3 is an extension of the ascent in Proposition 1.1(a). However the descent in Proposition 1.1(b) can not be extended in case where $\ell_R(S/\frak mS)>1$ (see Remark \ref{R:1a}).

In Section 2, we give some preliminaries that will be used in the sequel. In the last section, we prove Proposition 1.1 and Theorems 1.2, 1.3  (see Proposition \ref{P:1b} and Theorems \ref{T:1}, \ref{T:2}). As an application, we clarify the Artinianess of certain local cohomology modules (Corollaries \ref{C:1}, \ref{C:2}, \ref{C:3}).

\section{Preliminaries}

\ \ \ \ We recall some known properties of Artinian modules. Let $A$ be an Artinian $R$-module. Then $A$ is $\frak m$-torsion, i.e. $A=\bigcup_{t\in \Bbb N}(0:_A\frak m^t).$ Therefore, $A\neq 0$ if and only if $(0:_A\frak m)\neq 0$. 

A criterion of Artinian modules is given by L. Melkersson \cite[Theorem 1.3]{Mel}. 

\begin{lemma} \label{L:bx}  Let $A$ be an $R$-module. Then $A$ is Artinian if and only if $A$ is $I$-torsion and $(0:_AI)$ is Artinian for some ideal $I$ of $R$.
\end{lemma}

Let $A\neq 0$ be an $R$-module and $\frak \p\in\Spec(R).$ Following Macdonald \cite{Mac}, $A$ is said to be {\it $\frak p$-secondary} if  the multiplication by $a$ on $A$ is  surjective for all $a\in R\setminus \frak p$ and  it is nilpotent for all $a\in\frak p$.  We say that $A$ is {\it representable} if $A$ has a minimal secondary representation $A=A_1+\ldots +A_n$, where each $A_i$ is $\frak p_i$-secondary,  $A_i$ is not redundant and $\frak p_i\neq \frak p_j$ for all $i\neq j.$ The set $\{\frak p_1, \ldots , \frak p_n\}$ does not depend on the choice of the minimal secondary representation of $A$. This set is denoted by $\Att_R(A)$ and is called the set of {\it  attached primes} of $A$.

Note that any Artinian module is representable. Here are some elementary properties of attached primes of Artinian modules (see \cite{Mac}).  

\begin{lemma}\label{L:1} Let $A\neq 0$ be an Artinian $R$-module.  Then

(a) $\ell_R(A)<\infty$ if and only if $\Att_RA=\{\frak m\}.$ 

(b)  $\min\Att_RA=\min\Var (\Ann_RA).$ 
\end{lemma}

Denote by $\widehat R$ the $\frak m$-adic completion of $R$. We recall the ascent and descent of Artinian module structures between $R$ and $\widehat R$, see \cite[8.2.4, 10.2.18]{BS}.

\begin{remark} \label{R:1} {\rm Let $j: R\rightarrow \widehat R$ be the natural map.

(a) Let $A$ be an Artinian $R$-module. Then $A$  has a natural structure as an Artinian $\widehat R$-module, where each subset of $A$ is an $R$-submodule if and only if it is an $\widehat R$-submodule. If we regard this $\widehat R$-module as an $R$-module by means of $j$, then we recover the original $R$-module structure on $A$. Moreover, the map $A\rightarrow A\otimes_R\widehat R$ sending  $a$ to  $a\otimes 1$ is an isomorphism of $\widehat R$-modules. 

(b) Let $B$ be an Artinian $\widehat R$-module. Then $B$ is Artinian as an $R$-module by means of $j$. If we regard the Artinian $R$-module $B$ as an $\widehat R$-module using the method in (a), we recover the original $\widehat R$-module structure on $B$.}
\end{remark}

  It is clear that if $B$ is a $S$-module such that $\ell_R(B)<\infty$, then $\ell_S(B)<\infty$.  The following lemma gives the descent of finite length module structures between $R$ and $S$. 

\begin{lemma} \label{L:1a} The following statements are equivalent:

(a)  Every  finite length $S$-module is of finite length as an $R$-module (by means of $\varphi$).

(b) Every  finite length $S$-module is Artinian as an $R$-module.

(c) There exists a finite length $S$-module $B\neq 0$  such that $B$ is Artinian as an $R$-module.

(d) $\ell_R (S/\frak n)<\infty$.

Suppose that the equivalent conditions (a), (b), (c), (d) satisfy. If $B$ is a finite length $S$-module, then 
$$\ell_R(B)=\ell_S(B)\ell_R(S/\frak n).$$
\end{lemma}

\begin{proof} (a)$\Rightarrow$(b) and (b)$\Rightarrow$(c) are clear.

 (c)$\Rightarrow$(d).  Let $B\neq 0$ be a finite length $S$-module such that $B$ is an Artinian $R$-module. It is clear that $\dim (R/\Ann_R(0:_B\frak m))=0.$ It follows by Lemma \ref{L:1} that $\ell_R(0:_B\frak m)<\infty$. Therefore,  $\ell_S(0:_B\frak mS)<\infty$. Since $B\neq 0$, we have $(0:_B\frak mS)\neq 0$. So, $\Ass_S(0:_B\frak mS)=\{\frak n\}.$  So, $S/\frak n$ is isomorphic to a submodule of $(0:_B\frak mS).$ Therefore, $$\ell_R(S/\frak n)\leq \ell_R(0:_B\frak mS)<\infty.$$ 
\ \ \ (d)$\Rightarrow$(a).  Let $B$ be a finite length $S$-module.  The case $B=0$ is clear. Let  $B\neq 0$. Set  $\ell_S(B)=n.$ Then $n>0$  and there exists a chain $0=B_0\subset B_1\subset \ldots \subset B_n=B$ of $S$-submodules of $B$ such that  $B_i/B_{i-1}\cong S/\frak n$ for all $i=1, \ldots, n.$ By assumption (d) we get
$$\ell_R(B)=\sum_{i=1}^n \ell_R(B_i/B_{i-1})=n.\ell_R(S/\frak n)<\infty.$$
\end{proof}   

Since $\varphi$ is flat local, $S$ is a  faithfully flat $R$-module. Therefore, an $R$-module $M$ is finitely generated if and only if $M\otimes_RS$ is a finitely generated $S$-module.  The following proposition shows that the analogy for Artinian modules does not hold true any more.  Note that, for any $R$-module $A$, if $A\otimes_RS$ is  Artinian as a $S$-module, then $A$ is  Artinian as an $R$-module.

\begin{proposition} \label{P:1a} The following statements are equivalent:

(a)  $\dim(S/\frak mS)=0.$

(b)   $A\otimes_RS$ is an Artinian $S$-module for every Artinian $R$-module $A$. 

(c) $A\otimes_RS$ is an Artinian $S$-module for some Artinian $R$-module $A\neq 0$.

(d) $A\otimes_RS$ is a finite length $S$-module for every finite length $R$-module $A$.
 
(e) $A\otimes_RS$ is a finite length $S$-module for some finite length $R$-module $A\neq 0$.
\end{proposition}

\begin{proof} (a)$\Rightarrow$(b) follows by \cite[Lemma 2.3]{NQ}. The three assertions (b)$\Rightarrow$(c), (d)$\Rightarrow$(e) and (e)$\Rightarrow$(c) are clear.

(c)$\Rightarrow$(a). Let $A\neq 0$ be an Artinian $R$-module such that $A\otimes_RS$ is an Artinian $S$-module. Since $A$ is Artinian, $A$ is the direct limit of a direct system $\{A_n\}$ where each $A_n$ is a finite length $R$-submodule of $A$ and $A_n\neq 0$ for $n\gg 0$.  Hence $\Ass_R(A_n)=\{\frak m\}$ for $n\gg 0$. Note that tensor product commutes with direct limit, therefore we have by \cite[Theorem 23.2]{Mat} that
 $$\Ass_S(A\otimes_RS)\subseteq \underset{n}{\bigcup}\Ass_S(A_n\otimes_RS)=\Ass (S/\frak mS).$$
Since $S$ is a faithfully flat $R$-module, we have $\Ass_S(A_n\otimes_RS)\subseteq \Ass_S(A\otimes_RS)$ for all $n$. Hence $\Ass_S(A\otimes_RS)=\Ass (S/\frak mS).$ Since $A\otimes_RS\neq 0$ is an Artinian $S$-module, $\Ass_S(A\otimes_RS)=\{\frak n\}.$ Hence $\dim (S/\frak mS)=0.$

(a)$\Rightarrow$(d). Let $A$ be a finite length $R$-module. Then we have by the faithful flatness of $S$ over $R$ and by our assumption (a) that $\ell_S(A\otimes_RS)=\ell_R(A)\ell_S(S/\frak mS)<\infty$.
\end{proof}

\section{Main results}

\ \ \ \ Firstly, we give characterizations for the descent of Artinian module structures between $R$ and $S$.   Let $E(S/\frak n)$ be the injective hull of $S/\frak n$ over $S$. Let $r=\dim(S/\frak mS)$ and $H^r_{\frak n}(S/\frak mS)$ be the top local cohomology module of $S/\frak mS$ with respect to $\frak n$. Note that   $E(S/\frak n)$ and $H^r_{\frak n}(S/\frak mS)$ are Artinian $S$-modules, see  \cite[Theorems 7.1.3, 10.2.5]{BS}.

\begin{theorem}  \label{T:1} The following statements are equivalent:

(a)  Each Artinian $S$-module is Artinian as an $R$-module (by means of $\varphi$).

(b) $E(S/\frak n)$ is Artinian as an $R$-module.

(c) $H^r_{\frak n}(S/\frak mS)$ is Artinian as an $R$-module, where $r=\dim (S/\frak mS)$.

(d) $\ell_R(S/\frak m S)<\infty$.
\end{theorem}

\begin{proof} (a)$\Rightarrow$(b)  and (a)$\Rightarrow$(c) are clear.

(b)$\Rightarrow$(d).  Set $B=E(S/\frak n)$.  By assumption (b), $B$ is  non-zero Artinian as an $R$-module. Hence $0<\ell_R(0:_B\frak m)<\infty$ by Lemma \ref{L:1}.  Therefore, $0<\ell_S(0:_B\frak mS)<\infty$. Hence $\Rad (\Ann_S(0:_B\frak mS))=\frak n.$ We claim that $P= \Ann_S(0:_BP)$ for all  $P\in\Spec(S)$.  In fact, it is clear that if $P=\frak n$, then $P= \Ann_S(0:_BP)$. Let $P\in\Spec(S)$ such that $P\neq \frak n$. Then 
$$\Ass_S (\Hom_S(S_P, B))=\{Q\in\Spec(S)\mid Q\subseteq P\}$$
by \cite[Lemma 4.1]{MS} (see also \cite[Folgerlung 4.7]{Zos}). Hence $P\in \Ass_S (\Hom_S(S_P, B)).$ Let $0\neq f\in \Hom_S(S_P, B)$ such that $P=\Ann_S(f)$. Then we get $Pf(x)=0\subset B$, and hence $f(x)\in (0:_BP)$ for all $x\in S_P.$ Therefore, we have the induced homomorphism of $S$-modules  $f*: S_P\rightarrow  (0:_BP))$ which is defined by $f^*(x)=f(x)$ for all $x\in S_P.$ Hence $\Hom_S(S_P, (0:_BP))\neq 0$. By \cite[Lemma 7.2(2)]{MS}, there exists $Q\in \Att_S(0:_BP)$ such that $Q\subseteq P.$  Note that $\min\Att_S(0:_BP)=\min\Var (\Ann_S(0:_BP))$ by Lemma \ref{L:1}(b), it follows  that $P\supseteq \Ann_S(0:_BP).$ Therefore, $P= \Ann_S(0:_BP)$ and the claim is proved.  Now, let $P\in\Spec(S)$ such that $P\supseteq\frak mS.$ Then   $\Ann_S(0:_B\frak mS)\subseteq \Ann_S(0:_BP)=P$  by the claim. Hence $$\frak n=\Rad (\Ann_S(0:_B\frak mS))\subseteq \underset{P\supseteq \frak mS}{\underset{P\in\Spec(S)}\bigcap} P=\Rad (\frak m S).$$ Therefore, $\dim (S/\frak mS)=0.$ Since the finite length $S$-module $(0:_B\frak mS)\neq 0$ is an Artinian $R$-module, we get by Lemma \ref{L:1a}(c)$\Rightarrow$(d) that $\ell_R(S/\frak n)<\infty$. It follows again by Lemma \ref{L:1a} that $\ell_R(S/\frak mS)=\ell_S(S/\frak m S)\ell_R(S/\frak n)<\infty$.

(d)$\Rightarrow$ (a). Let $B$ is an Artinian $S$-module. We use Melkersson's criterion (see Lemma \ref{L:bx})  to prove that $B$ is Artinian as an $R$-module. Let $b\in B$. Since $B$ is an Artinian $S$-module, $\frak n^tb=0$ for some $t\in\Bbb N$.  Hence $\frak m^tb=0.$ Therefore,  $B$ is an $\frak m$-torsion $R$-module.  We get by assumption (d) that $\ell_S (S/\frak mS)<\infty$ and $\ell_R(S/\frak n)<\infty$. So, $\dim (S/\frak mS)=0$. Hence $\dim (S/\Ann_S(0:_B\frak mS))=0$ and hence $\ell_S(0:_B\frak mS)<\infty$ by Lemma \ref{L:1}. So, we have by Lemma \ref{L:1a}(d)$\Rightarrow$(b) that $(0:_B\frak m)$ is an Artinian $R$-module. Hence $B$ is Artinian as an $R$-module by Melkersson's criterion.

(c)$\Rightarrow$ (d). Set $B=H^r_{\frak n}(S/\frak mS)$, where $r=\dim(S/\frak mS)$. By our assumption (c), $B$ is non-zero Artinian as an $R$-module. It is clear that $\frak m=\Ann_R(B)$. Therefore, $\ell_R(B)<\infty$ by Lemma \ref{L:1} and hence $\ell_S(B)<\infty.$ So, $\ell_R(S/\frak n)<\infty$ by Lemma \ref{L:1a}(c)$\Rightarrow$(d) and $B$ is a finitely generated $S$-module. Therefore, $\dim (S/\frak mS)=0$ by \cite[Corollary 7.3.3]{BS}). Hence $\ell_R(S/\frak m S)<\infty$ by Lemma \ref{L:1a}.
\end{proof}

The following example shows that in the statements (b), (c) of Theorem \ref{T:1}, the Artinian $S$-module $E(S/\frak n)$ or $H^r_{\frak n}(S/\frak mS)$ can not be replaced by an arbitrary Artinian $S$-module. 

\begin{example} \label{E:1} {\rm Let $(R,\frak m)$ be a Noetherian local ring. Let $S=R[[x]]$ be the formal power series ring.  Let $\varphi: R\rightarrow S$ be the natural injection. Then $\varphi$ is a flat local homomorphism. Note that $\dim (S/\frak mS)=1$. So, $S/\frak mS$ is not of finite length as a $S$-module, and hence it is not of finite length as an $R$-module. Let $\frak n$ be the unique maximal ideal of $S$, set $t:=\depth (R).$ It is clear that $H^t_{\frak n}(S/xS)$ is an Artinian $S$-module. Note that $H^t_{\frak n}(S/xS)\cong H^t_{\frak m}(R)$ as $R$-modules. Therefore, $H^t_{\frak n}(S/xS)\neq 0$ and it is  Artinian as an $R$-module.}  
\end{example}

The following example illustrates the results in  Proposition \ref{P:1a}, Theorem \ref{T:1}.

\begin{example} \label{E:1a} {\rm  Let $R=k[[x_1, \ldots , x_d]]$ and $S=k'[[x_1, \ldots , x_{d'}]]$ be  the formal power series rings over  $k$ and  $k'$ respectively, where  $k$ is a field, $k'$ is a field extension of $k$ and $d'\geq d$. Then $\frak m=(x_1, \ldots , x_d)R$ and $\frak n=(x_1, \ldots , x_{d'})S$ are respectively the maximal ideals of $R$ and $S$.  It is clear that the natural injection $R\rightarrow S$ is a flat local homomorphism. We consider two cases:

(a) {\it The case $d=d'$}: Then $\dim(S/\frak mS)=0$. So, $A\otimes_RS$ is an Artinian $S$-module for any Artinian $R$-module $A$ (by Proposition \ref{P:1a}). If $k'$ is a finite extension of $k$ (for instance, $k=\Bbb Q$ and $k'=\Bbb Q[\alpha]$ with $\alpha$ an algebraic number) then $\ell_R(S/\frak mS)=\dim_k(k')<\infty$ and therefore each Artinian $S$-module is Artinian as an $R$-module (by Theorem \ref{T:1}). If $k'$ is an infinite extension of $k$ (for instance, $k=\Bbb Q$ and $k'=\Bbb Q(\alpha)$ with $\alpha$ a transcendental number) then $\ell_R(S/\frak n)=\infty$, therefore any non zero $S$-module is of infinite length as an $R$-module (by Lemma \ref{L:1a}) and $E(S/\frak n)$ is Artinian as a $S$-module but it is not Artinian as an $R$-module (by Theorem \ref{T:1}).

(b) {\it The case $d'>d$}: Then $\dim (S/\frak mS)>0$. So, $A\otimes_RS$ is not an Artinian $S$-module for any $R$-module $A\neq 0$ (by Proposition \ref{P:1a}). Moreover, $E(S/\frak n)$ and $H^{d'-d}_{\frak n}(S/\frak mS)$ are Artinian $S$-modules, but they are not Artinian as  $R$-modules (by Theorem \ref{T:1}). Note that $k'$ is a finite extension of $k$ if and only if $\ell_R(S/\frak n)<\infty$, if and only if any finite length  $S$-module is a finite length $R$-module (by Lemma \ref{L:1a}).}
\end{example}

As an application,  we study the Artinianess of the top local cohomology modules.

\begin{corollary} \label{C:1} Suppose that $S$ is a catenary Noetherian local domain of dimension $d$. Then $H^d_{\frak n}(S)$ is an Artinian $R$-module  by means of $\varphi$ if and only if $\ell_R(S/\frak mS)<\infty$.
\end{corollary}

\begin{proof} Set $B:=H^d_{\frak n}(S).$ Then $B\neq 0$ is an Artinian $S$-module and $\Att_S(B)=\{0\}$ by \cite[Theorem 7.3.2]{BS}. So, it follows by Lemma \ref{L:1}(b) that $\Ann_S(B)=0.$  If $\ell_R(S/\frak mS)<\infty$, then $B$ is  Artinian as an $R$-module by Theorem \ref{T:1}. Conversely, suppose that $B$ is Artinian as an $R$-module. Then $0<\ell_R(0:_B\frak m)<\infty$ by Lemma \ref{L:1}, hence $0<\ell_S(0:_B\frak mS)<\infty$. Therefore, $\frak n=\Rad(\Ann_S(0:_B\frak mS))$ and $\ell_R(S/\frak n)<\infty$ by Lemma \ref{L:1a}(c)$\Rightarrow$(d).  For any $P\in\Spec(S),$ since $P\supseteq 0=\Ann_S(B)$ and $S$ is catenary, we get by \cite[Theorem 4.1]{CDN} that $\Ann_S(0:_BP)=P.$ Hence 
$$\frak n=\Rad(\Ann_S(0:_B\frak mS))\subseteq \underset{P\supseteq \frak mS}{\underset{P\in\Spec(S)}{\bigcap}}\Rad(\Ann_S(0:_BP))=\underset{P\supseteq \frak mS}{\underset{P\in\Spec(S)}{\bigcap}}P=\Rad (\frak mS).$$ So, $\dim (S/\frak mS)=0.$ Therefore, $\ell_R(S/\frak mS)=\ell_S(S/\frak mS)\ell_R(S/\frak n)<\infty$ by Lemma \ref{L:1a}.
\end{proof}

 The {\it cohomological dimension} of an $R$-module $M$ with respect to an ideal $I$ of $R$, denoted by ${\rm cd}\ (I, M)$, is the largest integer $i$ such that $H^i_I(M)\neq 0$. 

\begin{corollary}\label{C:2}  Let $M$ be a finitely generated $R$-module of dimension $d$. Then  we have ${\rm cd}(\frak mS, M\otimes_RS)=d$ and  the top local cohomology module $H^d_{\frak mS}(M\otimes_RS)$ is  an Artinian $S$-module if and only if $\dim (S/\frak mS)=0$. 
\end{corollary}

\begin{proof} By the flat base change  \cite[Theorem 4.3.2]{BS}, we have $H^i_{\frak mS}(M\otimes_RS)\cong H^i_{\frak m}(M)\otimes_RS$ for all $i$. Hence ${\rm cd}(\frak mS, M\otimes_RS)=d$ by the faithful flatness of $S$ over $R$. The rest assertion follows by Proposition \ref{P:1a}.  
\end{proof}

 Next, we claim that if $\ell_R(S/\frak mS)=1$ then the ascent and descent of Artinian module structures between $R$ and $S$ are positive.  Note that  if $S=\widehat R$ or $S$ is the Henselization $R^h$ of $R$ introduced by M. Nagata \cite{Na} with the natural maps $\varphi$, then $\ell_R(S/\frak mS)=1$.

\begin{proposition}  \label{P:1b} Assume that $\ell_R(S/\frak mS)=1$. Let $A$ be an Artinian $R$-module, let $B$ be an Artinian $S$-module. The following statements are true. 

(a) $A$ has a natural structure as an Artinian $S$-module which is compatible with its original $R$-module structure, the map $A\rightarrow A\otimes_RS$ sending $a$ to $a\otimes 1$ and the map $\Hom_R(S, A)\to A$ sending $f$ to $f(1)$ are isomorphisms of $S$-modules, and $\Ext^i_R(S, A)=0$ for all integers $i>0.$

(b) $B$ is Artinian as an $R$-module by means of $\varphi$.  If we regard this $R$-module as a $S$-module using the method in (a), then  we recover the original $S$-module structure on $B$. 
\end{proposition}

\begin{proof} (a) Since $S$ is faithfully flat as an $R$-module,  the map $\varphi_t: R/\frak m^t\rightarrow S/\frak m^tS$ given by $\varphi_t(r+\frak m^t)=r+\frak m^tS$  is an injection for all positive integers $t$. As $\ell_R(S/\frak mS)=1$, we get $\frak n=\frak mS$. Hence $\ell_R(S/\frak n)=1$ and $\ell_S(S/\frak mS)=1.$ So,  $\ell_R(S/\frak m^tS)=\ell_S(S/\frak m^tS)$ by Lemma \ref{L:1a}. Therefore, for any $t\in\Bbb N$ we have  by the faithful flatness of $S$ over $R$ that
$$\ell_R(S/\frak m^tS)=\ell_S(R/\frak m^t\otimes_RS)=\ell_R(R/\frak m^t)\ell_S(S/\frak mS)=\ell_R(R/\frak m^t).$$ Therefore, $\varphi_t$ is an isomorphism. Let $s\in S$ and $a\in A.$ Then there exists $t\in\Bbb N$ such that $\frak m^ta=0$.  Since $\varphi_t$ is isomorphic, there exists uniquely an element $r+\frak m^t\in R/\frak m^t$ such that $\varphi_t(r+\frak m^t)=s+\frak m^tS.$ Set $s\circ a:=ra$.  Note that the definition of $s\circ a$ does not depend on the choice of $r$ and $t$. Therefore, it gives a scalar multiplication on $A$ so that $A$ is a $S$-module. With this structure, each subset of $A$ is a $S$-submodule if and only if it is an $R$-submodule. Therefore, $A$ is Artinian as a $S$-module. Moreover, if we regard this $S$-module as an $R$-module by means of  $ \varphi$, then we recover the original $R$-module structure on $A$.

Since $S$ is faithfully flat over $R$, the homomorphism $g: A\rightarrow A\otimes_RS$ sending $a$ to $a\otimes 1$ is an injection of $R$-modules. Consider the bilinear map $f: A\times S\rightarrow A$ given by $f(a,s)=s\circ a.$ Then there exists a homomorphism $h: A \otimes_RS\rightarrow A$ of $R$-modules given by $h(a\otimes s)=s\circ a.$ Let $a\in A, s\in S$. Let $t\in\Bbb N$ such that $\frak m^ta=0.$ Then $\frak m^tS(a\otimes 1)=0.$ Take $r\in R$ such that $\varphi_t(r+\frak m^t)=s+\frak m^tS.$ By the definition of the scalar multiplication of $A$ over $S$, we have 
$$gh(a\otimes s)=g(s\circ a)=g(ra)=ra\otimes 1=r(a\otimes 1)=s(a\otimes 1)=a\otimes s.$$
Hence $gh={\rm id}_{A\otimes_RS}.$ Therefore, $g$ is surjective and hence $g$ is an isomorphism of $R$-modules. We can check that $g$ is also an isomorphism of $S$-modules.

 Consider the $R$-homomorphism $p:\Hom_R(S,A)\to A$ sending $f$ to $f(1)$. For each $a\in A$, we define the map $f_a: S\rightarrow A$ given by $f_a(s)=s\circ a.$ Then $f_a\in \Hom_R(S,A).$ Then the map $q: A\to\Hom_R(S,A)$ sending $a$ to the map $f_a$ is a homomorphism of $R$-modules. We have $pq(a)=p(f_a)=f_a(1)=1\circ a=a$ for all $a\in A.$ Therefore $pq={\rm id}_A.$ Let $f\in\Hom_R(S, A)$ and $s\in S.$ Take $t\in\Bbb N$ and $r\in R$ such that $\frak m^tf(1)=0$ and $\varphi_t(r+\frak m^t)=s+\frak m^tS.$ Hence $f(\frak m^tS)=0$ and hence $f(r)=f(s).$ Therefore,
 $$(qp(f))(s)=(q(f(1))(s)=s\circ f(1)=rf(1)=f(r)=f(s).$$ It follows that $qp={\rm id}_{\Hom_R(S,A)}.$ So, $p$ is an isomorphism of $R$-module.   We can check that $p$ is also an isomorphism of $S$-modules.

  Let $ 0\to E_0\to E_1\to\ldots$ be an injective resolution of $A$ such that $E_i$ is Artinian $R$-module for all $i$. Since $\Hom_R(S,E_i)\cong E_i$ for all $i$, we obtain  $\Ext^i_R(S, A)=0$ for all $i>0.$
 
(b) Since $\ell_R(S/\frak mS)<\infty$, we get by Theorem \ref{T:1} that $B$ is Artinian as an $R$-module. The rest statement follows obviously. 
\end{proof}

  The following theorem is the second main result of this paper which clarifies the structure of $A\otimes_RS$ for each Artinian $R$-module $A$ in case where $\ell_R(S/\frak mS)=m<\infty$. Note that  $A^m$ is an Artinian $R$-module with the usual scalar multiplication $r(a_1, \ldots , a_m):=(ra_1, \ldots , ra_m)$ for all $r\in R$ and all $(a_1, \ldots , a_m)\in A^m.$   This theorem also shows that each Artinian $R$-module structure on $A$ induces an Artinian $S$-module structure on $A^m$ compatible with its usual $R$-module structure. 

\begin{theorem} \label{T:2} Assume that $\ell_R(S/\frak mS)=m$. Let $A$ be an Artinian $R$-module. Then  $A^m$ has a structure of Artinian $S$-module and there is an isomorphism $A\otimes_RS\cong A^m$ of $S$-modules. If we regard this $S$-module $A^m$ as an $R$-module by means of $\varphi$, then we recover the usual $R$-module structure on $A^m$.  
\end{theorem}

\begin{proof} Note that $S/\frak mS$ is a $R/\frak m$-vector space of dimension $m$.  Let $\epsilon_1, \ldots , \epsilon_m$ be elements in $S$ such that $\epsilon_1+\frak mS, \ldots , \epsilon_m+\frak mS$ is a  base of  the vector space $S/\frak mS$. We claim that $S/\frak m^tS=\big(R(\epsilon_1, \ldots , \epsilon_m)+\frak m^tS\big)/\frak m^tS$ for all $t\in\Bbb N$. Indeed, since $S=R(\epsilon_1, \ldots , \epsilon_m)+\frak mS,$ we have 
$$S/\frak m^tS=\big(R(\epsilon_1, \ldots , \epsilon_m)+\frak m^tS\big)/\frak m^tS + \frak m(S/\frak m^tS).$$ 
It follows by Lemma \ref{L:1a} and by the faithful flatness of $S$ over $R$ that 
\begin{align}\ell_R(S/\frak m^tS)=\ell_S(S/\frak m^tS)\ell_R(S/\frak n)&=\ell_R(R/\frak m^t)\ell_S(S/\frak mS)\ell_R(S/\frak n)\notag\\
&=\ell_R(R/\frak m^t)\ell_R(S/\frak mS)=m\ell_R(R/\frak m^t)=\ell_R\big((R/\frak m^t)^m\big)<\infty.\notag\end{align} Therefore, the claim follows by Nakayama Lemma. We consider  $S/\frak m^tS$ as an $R/\frak m^t$-module.  Then, by the claim we can define a surjection  $\varphi_t: (R/\frak m^t)^m\rightarrow S/\frak m^tS$ given by $$\varphi_t(r_1+\frak m^t, \ldots , r_m+\frak m^t)=(r_1\epsilon_1+ \ldots +r_m\epsilon_m)+\frak m^tS$$ for all $r_1, \ldots , r_m\in R.$ Since $\ell_R(S/\frak m^tS)=\ell_R\big((R/\frak m^t)^m\big),$ it follows that $\varphi_t$ is an isomorphism.

  Consider $A^m$ as an $R$-module with the usual scalar multiplication. Define a bilinear map $f: A\times S\rightarrow A^m$ as follows: Let $a\in A$ and $s\in S.$ Since $A$ is Artinian, there exists $t\in\Bbb N$ such that $\frak m^ta=0$. By the claim, there exist $r_1, \ldots , r_m\in R$ such that $s=\sum_{i=1}^m r_i\epsilon_i \ ({\rm mod} (\frak m^tS))$. Set $f(a, s):=(r_1a, \ldots , r_ma).$ As $\varphi_t$ is isomorphic,  we can check that the  definition of  $f(a, s)$ does not depend on the choice of  $t$ and $r_1, \ldots , r_m$.  Therefore, there exists an $R$-homomorphism $h: A \otimes_RS\rightarrow A^m$ given by $h(a\otimes s)=f(a, s)$ for all $a\in A$ and $s\in S.$ Let $(a_1, \ldots , a_m)\in A^m$, we have 
 $$(a_1, \ldots , a_m)=\sum_{i=1}^m\big(f(a_i,\epsilon_i)\big)=\sum_{i=1}^m h(a_i\otimes \epsilon_i)=h\big(\sum_{i=1}^m(a_i\otimes \epsilon_i)\big).$$ Therefore, $h$ is a surjection. 

We consider a $R$-homomorphism $g: A^m\rightarrow A\otimes_RS$  given by $g(a_1, \ldots , a_m):=\sum_{i=1}^m(a_i\otimes\epsilon_i)$ for any $(a_1, \ldots , a_m)\in A^m$. Let $a\in A, s\in S$ and let  $t\in\Bbb N$ be such that $\frak m^ta=0$.  By the claim, there exist $r_1, \ldots , r_m\in R$ such that $s=\sum_{i=1}^m r_i\epsilon_i\ ({\rm mod}(\frak m^tS)).$ Note that 
$a\otimes s=a\otimes \sum_{i=1}^m (r_i\epsilon_i)=\sum_{i=1}^m r_i(a\otimes \epsilon_i).$ Therefore,  $\{a\otimes \epsilon_i\mid a\in A, 1\leq i\leq m\}$ is a system of generators of $R$-module $A\otimes_R S$. By the definitions of $f, g, h$ we have 
$$gh(a\otimes \epsilon_1)=gf(a,\epsilon_1)=g(a,0, \ldots ,0)=a\otimes\epsilon_1.$$ Similarly, we have $gh(a\otimes \epsilon_i)=a\otimes\epsilon_i$ for all $2\leq i\leq m.$ It follows that $gh={\rm id}_{A\otimes_RS}.$ Therefore, $h$ is injective. Hence, $h$ is an isomorphism of $R$-modules and $hg={\rm id}_{A^m}.$ 

Now we provide a  $S$-module structure on $A^m$ as follows. Let $\underline a=(a_1, \ldots , a_m)\in A^m$ and $s\in S.$  Set $$s\circ\underline a:=h\big(s(g(\underline a))\big).$$
For any $s, s'\in S$, $\underline a, \underline a'\in A^m$, it is clear that $(s+s')\circ\underline a=s\circ\underline a+s'\circ\underline a$, $s\circ(\underline a+\underline a')=s\circ\underline a +s\circ\underline a'$. Since $hg={\rm id}_{A^m},$ we have $1\circ\underline a=h(1g(\underline a))=hg(\underline a)=\underline a.$ Moreover, 
$$s'\circ(s\circ\underline a)=h(s'g(s\circ\underline a))=h(s'g(h(sg(\underline a))))=h(s'sg(\underline a))=(s's)\circ\underline a.$$ So, this is a scalar multiplication and it makes $A^m$ being a $S$-module. Consider this $S$-module $A^m$ as an $R$-module by means of $\varphi$, then for $r\in R$ and $\underline a=(a_1, \ldots , a_m)\in A^m$ we have  
\begin{align}\varphi(r)\circ\underline a=h(\varphi(r)g(\underline a))&=h(\varphi(r)\sum_{i=1}^m(a_i\otimes\epsilon_i))=\sum_{i=1}^m h(a_i\otimes \varphi(r)\epsilon_i)\notag\\
&=\sum_{i=1}^m h(a_i\otimes r\epsilon_i)=\sum_{i=1}^m h(ra_i\otimes \epsilon_i)=(ra_1, \ldots , ra_m).\notag\end{align} Therefore, we recover the usual $R$-module structure on $A^m$. 

Let $s\in S$ and $z\in A\otimes_RS.$ Then $\displaystyle z=\sum_{i=1}^m u_i(a_i\otimes \epsilon_i)$ for some elements $u_i\in R$ and $a_i\in A$. By the definition of $h$,  $$h(z)=\sum_{i=1}^m h(u_ia_i\otimes \epsilon_i)=(u_1a_1, \ldots , u_ma_m).$$ Since $gh={\rm id}_{A\otimes_RS},$ we get by the definition of the scalar multiplication of $A^m$ over $S$ that 
 $$s\circ h(z)=s\circ (u_1a_1, \ldots , u_ma_m)=h(sg(u_1a_1, \ldots , u_ma_m))=h(sgh(z))=h(sz).$$ Therefore, $h$ is also an isomorphism of $S$-modules. Since $A\otimes_RS$ is an Artinian $S$-module by Proposition \ref{P:1a}, it follows that  $A^m$ is also an Artinian $S$-module.
\end{proof}

\begin{remark} \label{R:1a} {\rm Let  $\ell_R(S/\frak mS)=m$. Let $B$ be an Artinian $S$-module. Then $B^m$ is an Artinian $S$-module with the usual scalar multiplication $s(b_1, \ldots , b_m):=(sb_1, \ldots , sb_m).$  On the other hand, by Theorem \ref{T:1},  $B$ is Artinian as an $R$-module by means of $\varphi$. So, $B^m$ is an Artinian $S$-module with the second scalar multiplication $s\circ (b_1, \ldots , b_m):=h(sg(b_1, \ldots , b_m)),$ where $h, g$ are defined as in the proof of Theorem \ref{T:2}. If $m=1$, then clearly that the two scalar multiplications of $S$ on $B$ are the same. However, this is not the case for $m>1,$ i.e. the two $S$-module structures  on $B^m$ are different in general. For example, let $R=\Bbb Q$, $S=B=\Bbb Q[\sqrt 2]$ and $\varphi: R\rightarrow S$ is the natural injection. Then $m=2$ and we can choose $\epsilon_1=1$, $\epsilon_2=\sqrt 2$. For $s=1+3\sqrt 2\in S$ and $\underline b=(1, \sqrt 2)\in B^2$, the usual scalar multiplication gives $s\underline b= (1+3\sqrt 2, 6+\sqrt 2)$. Since $s\epsilon_1=1\epsilon_1+3\epsilon_2$ and $s\epsilon_2=6\epsilon_1+\epsilon_2$, the second scalar multiplication gives 
$$s\circ\underline b=h(sg(\underline b))=h(s(1\otimes\epsilon_1+\sqrt 2\otimes\epsilon_2))=h(1\otimes s\epsilon_1)+h(\sqrt 2\otimes s\epsilon_2)=(1+6\sqrt 2,3+\sqrt 2).$$}
\end{remark}

 The following result is a consequence of Theorem \ref{T:2}.
\begin{corollary}\label{C:3}  Let $M$ be a finitely generated $R$-module. If $\ell_R (S/\frak mS)=m$, then there is an isomorphism  $H^i_{\frak n}(M\otimes_RS)\cong (H^i_{\frak m}(M))^m$ of Artinian $S$-modules for all integers $i\geq 0$. 
\end{corollary}

\begin{proof} Note that $\Rad (\frak mS)=\frak n.$ Since $S$ is flat over $R$, it follows by Flat Base Change \cite[Theorem 4.3.2]{BS} and Theorem \ref{T:2} that 
$$H^i_{\frak n}(M\otimes_RS)= H^i_{\frak mS}(M\otimes_RS)\cong H^i_{\frak m}(M)\otimes_RS\cong (H^i_{\frak m}(M))^m$$ for all integers $i\geq 0$.
\end{proof}

Let $A\neq 0$ be an Artinian $R$-module.  Theorem \ref{T:2} describes the structure of $A\otimes_RS$ in case where $\ell_R(S/\frak mS)<\infty$. In the case where $\ell_R(S/\frak mS)=\infty$, if $\dim (S/\frak mS)=0$ then $A\otimes_RS$ is Artinian as a $S$-module by Proposition \ref{P:1a}, but $A\otimes_RS$ is not Artinian as an $R$-module by Theorem \ref{T:1}. So, it is natural to ask whether we can provide a $S$-module structure on $\bigoplus_{i\in I}A_i$ such that there is an isomorphism $A\otimes_RS\cong \bigoplus_{i\in I}A_i$ of $S$-modules,  where $I$ is a base of the $R/\frak m$-vector space $S/\frak mS$ and $A_i$ is a copy of $A$ for all $i\in I$. This question is still open now.

\end{document}